\documentclass[12pt]{amsart}
\usepackage{amssymb,amsmath,amsthm,slashed}
\newtheorem{theorem}{Theorem}
\newtheorem{lemma}[theorem]{Lemma}
\newtheorem{corollary}[theorem]{Corollary}

\newtheorem{proposition}[theorem]{Proposition}
\newtheorem{definition}[theorem]{Definition}
\newtheorem{fact}[theorem]{Fact}

\theoremstyle{remark}

\newtheorem{Remark}{Remark}
\newtheorem*{remark}{Remark}

\numberwithin{theorem}{section} \numberwithin{equation}{section}

\newcommand{\Z}{\mathbb{Z}}

\newcommand{\re}{\textnormal{Re}\,}
\newcommand{\im}{\textnormal{Im}\,}

\newcommand{\beq}{\begin{eqnarray*}}
\newcommand{\eeq}{\end{eqnarray*}}

\DeclareMathOperator{\DET}{DET}
\def\UPo{\overset{\circ}{\mathrm{UP}}}

\def\Zo{\overset{\circ}{\mathrm{Z}}}
\def\UP{\mathrm{UP}}
\def\Z{\mathrm{Z}}
\def\CP1{\mathbb{C}\mathrm{P}^1}
\def\Zc{\overline{\mathrm{Z}}}
\def\hoehe{$\phantom{\overset{\circ}{|}}$}

\title{The signature of the Seiberg-Witten surface}

\author{Andreas Malmendier}

\address{Department of Mathematics, Colby College, 
Waterville, ME 04901}
\email{andreas.malmendier@colby.edu}

\begin{document}

\begin{abstract}
The Seiberg-Witten family of elliptic curves defines a 
Jacobian rational elliptic surface $\Z$ over $\mathbb{C}\mathrm{P}^1$. We show that for the $\bar{\partial}$-operator along the fiber
the logarithm of the regularized determinant $- \frac{1}{2} \log \det' (\bar\partial^* \bar\partial)$ 
satisfies the anomaly equation of the one-loop topological string amplitude derived in Kodaira-Spencer theory.
We also show that not only the determinant line bundle with the Quillen metric but also the $\bar{\partial}$-operator itself 
extends across the nodal fibers of $\mathrm{Z}$. The extension introduces 
current contributions to the curvature of the determinant line bundle
at the points where the fibration develops nodal fibers.
The global anomaly of the determinant line bundle
then determines the signature of $\mathrm{Z}$ which equals minus the number of hypermultiplets.  
\end{abstract}

\maketitle

\section{Introduction and Statement of results}
There has been a continuing interest in  the non-perturbative properties of the supersymmetric Yang-Mills 
theory on four-dimensional manifolds. One of the results of the work of Seiberg and Witten \cite{SW2} 
is that the moduli space of the topological  $SU(2)$-Yang-Mills theory on a four-dimensional manifold 
decomposes into two branches, the Coulomb branch and the Seiberg-Witten branch. The branches are interpreted as 
the moduli spaces of simpler physical theories on the four-manifold. The Coulomb branch, also called the Seiberg-Witten family of curves, is the moduli 
space of a topological $U(1)$-gauge theory, called the low energy 
effective $U(1)$-gauge theory. We investigate the geometry and topology of the Coulomb branch as it is fundamental for the 
definition and the understanding of the $N=2$ supersymmetric, low energy effective field theory
by using the results and techniques developed in \cite{Atiyah2, SeeleyS, BismutB1, BismutB2}.

This article is structured as follows: In Section~\ref{geometry} we explain how
the Seiberg-Witten families of elliptic curves for the $N=2$ supersymmetric $SU(2)$-gauge theory with $N_f=0,1,2,3,4$ additional fields, called hypermultiplets, define four-dimensional, 
Jacobian rational elliptic surfaces $\Zc \to \CP1$ with singular fibers. We describe the correspondence
between the rational elliptic surfaces and the constraints on the masses of the hypermultiplets. 
Unless noted otherwise, we will then always assume that the masses are generic so that any singular fiber over $[u^*:1] \in \CP1$ for finite $u^*$ 
is a node and does not give rise to a surface singularity. Only the singular fiber over $[1:0]$ is a cusp giving rise to a surface singularity.
Thus, after removing the singular fiber over $[1:0]$, we obtain a smooth elliptic surface 
$\Z$ with three-dimensional boundary and elliptically fibered over a bounded disc in $\CP1$, called the $u$-plane $\UP \subset \CP1$.

In Section~\ref{single_curve} we review the construction of the regularized determinant of the $\bar{\partial}$-operator
on an elliptic curve. When the elliptic curve $E_u$ is varied in the Seiberg-Witten
family over the $u$-plane, we obtain the determinant line bundle $\DET \bar{\partial} \to \UP$ of the $\bar{\partial}$-operator along the fiber of $\Z \to \UP$. 
The regularized determinant $\sideset{}{^\prime}\det (\bar\partial^* \bar\partial)$ of the Laplacian along the fiber  
becomes a smooth function on $\UP$ which we will later use in the Quillen construction to define a metric and connection on $\DET \bar{\partial}$. 
In Section~\ref{string_amplitude}, we show that the logarithm of the regularized determinant of the Laplacian $-\frac{1}{2}\ln
\sideset{}{^\prime}\det (\bar\partial^* \bar\partial)$ satisfies the anomaly equation for the one-loop topological string amplitude of 
Kodaira-Spencer theory derived in \cite{BCOV}.

In Section~\ref{relative_determinant0} we show that the local anomaly of $\DET \bar{\partial}$ vanishes
and determine the non-trivial global anomaly as holonomy of the determinant section. 
We use the results of Bismut and Bost \cite{BismutB1, BismutB2} to 
show that the determinant line bundle with the Quillen metric extends smoothly across the nodal fibers.
Because of the non-trivial holonomy the extension of the determinant line bundle
introduces current contributions to the curvature over the points in the $u$-plane where the fiber develops a node.

In the case that the singular fibers are nodes Seeley and Singer \cite{SeeleyS} showed that there is a $\bar{\partial}$-operator that is defined
on the nodal fiber as well so that in a neighborhood $U \subset \UP$ the family of operators $\lbrace \bar{\partial}_u \rbrace_{u \in U}$ is a
continuous family. However, their procedure of obtaining the determinant line bundle is different from Bismut and Bost \cite{BismutB1, BismutB2} 
because their Laplacian is different. In Section~\ref{singular_dbar} we explain the connection between the two procedures by a local change in the conformal gauge 
of the fiber metric in a neighborhood of the nodal fiber.

In Section~\ref{relative_signature} and Section~\ref{signature_operator} we discuss the elliptic operators connected to the signature of 
the elliptic surface $\Z \to \UP$. First, for the signature operator $\mathrm{D}$ along the fiber of $\Z \to \UP$ 
we compute the global anomaly. We show that there is a canonical trivialization of the determinant line bundle $(\DET \; \mathrm{D})^{\otimes 6}$, and the well-defined logarithmic monodromies
of the canonical section of $\DET \; \mathrm{D}$ determine the signature of $\Z \to \UP$. We interpret the determinant line bundle as a solution to a Riemann-Hilbert problem on $\CP1$.
On the other hand, the generalization of Hirzebruch's signature theorem for manifolds with boundary by Atiyah, Patodi, Singer (APS) \cite{APS} 
shows that the elliptic signature operator on the four-dimensional surface $\Z$ has an 
analytic index if one imposes APS boundary conditions on $\partial \Z$. We show that this analytic index equals minus the number $N_f$ of hypermultiplets.

\section{The Jacobian elliptic surfaces for $N=2$ Yang-Mills theory} 
\label{geometry} 
An elliptic curve $E$ in the 
Weierstrass form can be written as
\begin{eqnarray}
 \label{jacobian} y^2 = 4 x^3 - g_2 \, x - g_3 \;, 
\end{eqnarray} 
where $g_2$ and $g_3$ are numbers such that the discriminant $\pmb{\Delta}
=g_2^3- 27 g_3$ does not vanish. In homogeneous coordinates $[X:Y:W]$
Eq.~(\ref{jacobian}) becomes 
\beq
  W Y^2 = 4 X^3 -  g_2 \, X W^2 - g_3 \, W^3.
\eeq
 One can check that the point $P$ with coordinates $[0:1:0]$ is a 
smooth point of the curve. We consider $P$ the base point of the elliptic 
curve and the origin of the group law on  $E$. The two types of singularities that can occur as 
Weierstrass cubic are a rational curve with a node, which appears 
when the discriminant vanishes and $g_2 ,g_3 \not = 0$, or a cusp when $g_2=g_3=0$.

Next, we look at a family of cubic curves over $\mathbb{C}\mathrm{P}^1$.
The family is parametrized by the base space $\CP1$ and a line bundle 
$\mathcal{N} \to \mathbb{C}\mathrm{P}^1$. The quantities 
$g_2$ and $g_3$ are promoted to global sections of 
$\mathcal{N}^{\otimes 4}$ and $\mathcal{N}^{\otimes 6}$ respectively; the 
discriminant becomes a section of $\mathcal{N}^{\otimes 12}$. If the
sections are generic enough so that they do not always lie in the 
discriminant locus, we obtain a Weierstrass fibration $\pi: \Zc \to \CP1$ with a section, called a Jacobian elliptic fibration.
Each fiber comes equipped with the base point $P$
that defines a section $S$ of the elliptic fibration which does not pass through the nodes or cusps. 
We will always assume that $\mathcal{N}=\mathcal{O}_{\mathbb{C}\mathrm{P}^1}(-1)$. 

In addition, we will assume from now on that in the coordinate chart $[u:1] \in \mathbb{C}\mathrm{P}^1$, the discriminant $\pmb{\Delta}$ is a polynomial of degree 
$2 \le N_f+2 \le 6$ in $u$  where $0 \le N_f \le 4$, and $g_2$ and $g_3$ are polynomials in $u$ of degree at most $2$ and $3$ respectively. The space of all such Weierstrass elliptic surfaces has $N_f+1$ moduli. 
To see this consider first the case where $N_f=4$. From the seven parameters defining $g_2$ and $g_3$ two can be eliminated by scaling and a shift in $u$. Furthermore, 
we can arrange the coefficient of $g_2$ of degree two and the coefficient of $g_3$ of degree three to be the modular invariants of an elliptic curve 
with periods $1$ and $\tau_0$. The remaining four coefficients can be expressed in terms of four complex parameters. 
Following the convention of \cite{SW2}
we will denote the parameters by $m_1, \dots, m_{N_f}$.  In physics, they are called the {\it masses of the hypermultiplets}.

A non-trivial elliptic fibration has to develop singular fibers.
The classification of the singular fibers is part of Kodaira's classification theorem of all possible singular fibers 
of an elliptic fibration \cite{Kodaira}. 
For generic values of the masses, the polynomial $\pmb{\Delta}$ has $N_f+2$ simple zeros $u^*_1, \dots, u^*_{N_f+2}$ for $|u|<\infty$
with $g_2(u^*_i), g_3(u^*_i) \not = 0$ for $i=1,..,N_f+2$. From Kodaira's classification theorem of singular fibers \cite{Kodaira}
it follows that the elliptic fibration develops the nodes, i.e., a singular fibers of Kodaira type $I_1$
over the points $u^*_1, u^*_2 , \dots, u^*_{N_f+2}$. 
For special values for $m_1, \dots, m_{N_f}$, several singular fibers of Kodaira type $I_1$ can coalesce and form singular fibers of Kodaira type $I_k$ with $k \ge 2$,
where the discriminant has a zero of order $k$. 
The second chart over the base space is $[1:v] \in \mathbb{C}\mathrm{P}^1$. The 
intersection of the two charts is given by $u=1/v$ with $v \not=0$.
The Weierstrass coordinates transform according to $x \mapsto v^2 x$ and $y \mapsto v^3 y$;
since $g_2$ and $g_3$ are sections of $\mathcal{N}^4$ and $\mathcal{N}^6$ respectively,
they transform according to $g_{2} \mapsto v^4 \, g_{2}$ and $g_{3} \mapsto v^6 \, g_{3}$.
The discriminant $\pmb{\Delta} \mapsto v^{12} \, \pmb{\Delta}$ becomes a polynomial in $v$ of degree $10-N_f$. 
From Kodaira's classification theorem it follows that the singular fiber $E_\infty$ over $u=\infty \; (v=0)$ is a cusp, a singular fiber of Kodaira type $I_{4-N_f}^*$. 
\begin{Remark}
\label{change_of_variables2}
Under the change of the coordinate chart from $[u:1]$ to $[1:v]$ on $\CP1$ by $u=-1/v$ with $v \not=0$
the holomorphic one-form $dz=dx/y$ transforms as
\beq
  dz_u = \frac{dx_u}{y_u} = - \, v \, \frac{dx_v}{y_v} = - \, v \, dz_v \;.
\eeq
It follows that $(dz)^{\otimes 2}$ has the same transformation under a coordinate change as $(du)^{-1}$
whence $\mathcal{N}^2 \cong \omega^{-1}_{\CP1}$ where
$\omega_{\CP1}$ is the canonical bundle on $\CP1$.
\end{Remark}

The elliptic surface $\Zc$ is a hyper-surface in the variables $(u,[X:Y:W])$. $\Zc$ has surface singular points whenever all partial derivatives in $u, x, y$ simultaneously vanish.
Singular fibers of Kodaira type $I_1$ do not give rise to surface 
singularities, whereas all singular fibers of Kodaira type $I_n$, 
with $n \ge 2$, and $I^*_{n}$, with $n\ge 0$, do. 
The monodromy around singular fibers of type $I_n$ or $I_n^*$ is parabolic \cite{Kodaira}.
It is known \cite[Sec.~4.6]{Miranda2} that a Weierstrass fibration is rational 
(i.e., birational to $\mathbb{C}\mathrm{P}^2$), 
if $g_2$ and $g_3$ are polynomials in $u$ of degree at most $4$ and $6$ respectively. 
The minimal resolution $\Z$ of $\Zc$ is the blow-up of $\mathbb{C}\mathrm{P}^2$ in nine points, and 
therefore has Picard number $10$. 
Conversely, by contracting every component of the fiber which
does not meet $S$, we obtain back the normal surface $\Zc$. 
The section  $S$ and smooth fiber use up two dimensions, and so the number of components of any singular fiber is at most eight since the components are 
always independent in the Neron-Severi group. However, as we will see below not all configurations of singular fibers exist. 

For later use we also introduce the following notation: we will denote  by $\UP$ the base curve $\CP1$ minus a small disc around $u=\infty$, and the restrition of the elliptic fibration to $\UP$ will be denoted by $\Z \to \UP$. Similarly, $\UPo$ will denote the base curve $\CP1$ with small open discs 
around all points with singular fibers removed, and the restriction of the elliptic fibration
to $\UPo$ will be denoted by $\Zo \to \UPo$.
  
\begin{Remark}
\label{canonical_bundle}
We denote by $\omega_\Z$ the canonical bundle of $\Z$. The space $H^0(\omega_\Z)$ is the space
of global holomorphic two-forms of dimension $p_g=h^{2,0}$. 
The bundle $\omega_{\Z/\CP1}= \omega_\Z \otimes (\pi^* \omega_{\CP1})^{-1}$ restricts to 
the canonical line bundle $K_u=\omega_{E_u}$ on each smooth fiber $E_u$.  It is well-know that for any rational elliptic surface we have $h^{1,0}=h^{2,0}=0$ \cite[Thm.~2.10]{Miranda2} and the canonical class is minus the fiber class.
Hence, we have $c^2_1(\Z)=0$.
It also follows that the first Chern class $c_1( \omega_{\Z/\CP1} )$ is a pullback from the base manifold and $c^2_1( \omega_{\Z/\CP1} )=0$.
\end{Remark}

\noindent
The above discussion motivates the following definition:
\begin{definition}
A Seiberg-Witten curve for $N_f$ hypermultiplets is a Jacobian rational elliptic surface with one singular fiber of Kodaira type $I^*_{4-N_f}$
and singular fibers of Kodaira type $I_n$ and $I_n^*$ only.
\end{definition}
Using the explicit Weierstrass parametrization given in \cite{SW2} for $0 \le N_f \le 4$ it is easy to compose a
list of the possible configurations of singular fibers that appear as Seiberg-Witten curves. 
In Table \ref{table1}, we list the constraints on the moduli, which substituted into the Weierstrass presentation 
in \cite{SW2} realize the configuration of singular fibers, the structure of the singular fibers $E_{u^*_n}$ over 
finite $u^*_n$, and $E_{\infty}$ for the rational elliptic surface $\Zc \to \CP1$. Here, $r=8-\sum_{\nu} (m_\nu-1)$ 
where the sum in $\nu$ runs over all singular fibers of the elliptic surface, and $m_{\nu}$ denotes the number
of irreducible components in the singular fiber.

\begin{remark}
The completeness of Table \ref{table1} follows by comparing the list with the list of impossible configurations in \cite{Miranda}. It follows from \cite{Nori} that among the Seiberg-Witten curves the ones with $r=0$ are the 
only modular elliptic surfaces.
\end{remark}

{\small
\begin{center}
\begin{table}[ht]
\caption{}
\begin{tabular}{|c|c|c|c|l|l|}
       \hline
       $N_f$ & $r$ & $E_\infty$ & singular fibers $E_{u^*_n}$ & mass constraints \\
       \hline
       4  & $4$ & $I_0^*$ & $6 I_1$        & - \\
       4  & $3$ & $I_0^*$ & $I_2, 4I_1$    &$m_3=m_4 \not = 0$ \\
       4  & $2$ & $I_0^*$ & $I_3, 3I_1$    &$m_2 = m_3 = m_4 \not = 0$ \\
       4  & $2$ & $I_0^*$ & $2I_2, 2I_1$   &$m_3=m_4=0$\\
       4  & $1$ & $I_0^*$ & $I_4, 2I_1$     &$m_2 = m_3 = m_4 =0 $ \\
       4  & $1$ & $I_0^*$ & $3I_2$             &$m_1 = m_2,  m_3 = m_4 =0$ \\
       4  & $0$ & $I_0^*$ & $I_0^*$            &$m_1 = m_2 = m_3 = m_4 = 0$\\
       \hline
       3  & $3$ & $I_1^*$ & $5 I_1$      & - \\
       3  & $2$ & $I_1^*$ & $I_2, 3I_1$  & $m_2=m_3$ \\
       3  & $1$ & $I_1^*$ & $I_3, 2I_1$  & $m_1 = m_2 = m_3 \not =0$\\
       3  & $1$ & $I_1^*$ & $2I_2, I_1$  & $m_1=m_2=0$ \\
       3  & $0$ & $I_1^*$ & $I_4, I_1$   & $m_1 = m_2 = m_3 = 0$\\
       \hline
       2  & $2$ & $I_2^*$ & $4I_1$       & - \\
       2  & $1$ & $I_2^*$ & $I_2, 2I_1$  & $m_1 = m_2 \not = 0$\\
       2  & $0$ & $I_2^*$ & $2I_2$       & $m_1 = m_2 = 0$\\
       \hline
       1  & $1$ & $I_3^*$ & $3I_1$ & - \\ 
       \hline
       0  & $0$ & $I_4^*$ & $2I_1$ & -\\
       \hline
      \end{tabular}
\label{table1}
\end{table}
\end{center}
}

\section{The regularized determinant on an elliptic curve}
\label{single_curve}
\label{determinant_on_single_curve}
We consider an elliptic curve $E$ with periods $2\, \pmb{\omega}$ and $2\, \pmb{\omega'}$, modular parameter  $\tau = \frac{\pmb{\omega'}}{\pmb{\omega}}$, and
complex coordinate $z$. Let $\xi = \xi^1 + i \xi^2$ be the complex coordinate on the normalized torus with periods $1$ and $\tau$ such that
$\xi= \frac{z}{2\,\pmb{\omega}}$. For $n_1, n_2 \in \mathbb{N}$, a complex function $\varphi$ on the normalized torus with the periodicities
\begin{eqnarray*}
 \varphi(\xi^1 + 1 , \xi^2) & = & - e^{\pi i \nu_1} \; \varphi(\xi^1,\xi^2) \;,\\
 \varphi(\xi^1 + \re \tau ,\xi^2 + \im \tau) & = & - e^{\pi i \nu_2} \; \varphi(\xi^1,\xi^2) \;,
\end{eqnarray*}
is given by
\begin{equation*}
\label{eigenfunctions}
\begin{split}
& \varphi_{n_1,n_2}(\xi^1,\xi^2) \\  
=  & \exp 2\pi i \left\lbrace \left[ n_1 + \frac{1-\nu_1}{2}\right] \xi^1 + \frac{1}{\im\tau}
 \left[ n_2 + \frac{1-\nu_2}{2} - \re\tau \; \left( n_1 + \frac{1-\nu_1}{2}\right) \right] \xi^2 \right \rbrace \;.
\end{split}
\end{equation*}
In fact, the set of functions $\{\varphi_{n_1,n_2}\}$ constitutes a complete system of eigenfunctions for the Laplace operator $-4\partial_\xi \bar\partial_\xi$ 
where $2\bar\partial_\xi = \partial_{\xi^1} + i \partial_{\xi^2}$. Their eigenvalues under $2\bar\partial_\xi$ are
\begin{equation}
\label{eigenvalues_dbar}
 \frac{2\pi}{\im\tau} \left\lbrace \left( n_1 + \frac{1-\nu_1}{2}\right)\tau - \left( n_2 + \frac{1-\nu_2}{2}\right)
 \right\rbrace \;.
\end{equation}
Because of $2\bar\partial = 2\bar\partial_z = \frac{1}{\overline{\omega}} \bar\partial_\xi$
the functions $\varphi_{n_1,n_2}$ are also eigenfunctions of $2\bar\partial$ for the eigenvalues
\[
 \frac{\pi}{ \im\tau \; \overline{\pmb{\omega}}} \left\lbrace \left( n_1 + \frac{1-\nu_1}{2}\right) \tau - \left( n_2 + \frac{1-\nu_2}{2}\right)
 \right\rbrace \;.
\]
The holomorphic line bundle of positive spinors on an elliptic curve $E$
can also be interpreted as a holomorphic square root $K^{1/2}$ of the bundle of holomorphic 
$(1,0)$-forms $K = \Omega^{1,0}(E)$. The chiral Dirac operators 
are
\beq
\slashed{\partial}^{+} = 
\bar{\partial}: \; C^\infty(K^{1/2}) & \to & C^\infty(K^{1/2}\otimes \overline{K}) \;,\\
 \slashed{\partial}^{-} = - \partial: \; C^\infty(\overline{K^{1/2}}) & \to & C^\infty(\overline{K^{1/2}} \otimes K) \;.
\eeq
Equivalently, we can view the situation as follows: we choose the unique even spin structure on $E$ as a reference square root $K_0^{1/2}$.
$K_0^{1/2}$ is the preferred spin bundle for the chosen homology basis, its divisor $\kappa = 1/2 + \tau/2$ is the vector of Riemann
constants. If we twist the Dirac operator $\slashed{\partial}^{+}$ by a flat holomorphic line bundle $W_{(\nu_1,\nu_2)}$ of order two with divisor 
$\check{\xi} = - \frac{\nu_2}{2} - \frac{\nu_1}{2} \tau$ and $\nu_1, \nu_2 \in \lbrace 0,1 \rbrace$,
the twisted chiral Dirac operator becomes
\begin{equation}
 \label{chiral_dirac}
 \slashed{\partial}^{+}_{(\nu_1,\nu_2)} = \bar{\partial}_{(\nu_1,\nu_2)}: 
 \; C^\infty \left( K_0^{1/2}  \otimes W_{(\nu_1,\nu_2)} \right)  \to  
    C^\infty \left( K_0^{1/2}  \otimes W_{(\nu_1,\nu_2)} \otimes \overline{K}\right) \;.
\end{equation}
In other words, the twisted chiral Dirac operator is the $\bar{\partial}$-operator coupled
to the holomorphic line bundle $L = K_0^{1/2} \otimes W_{(\nu_1, \nu_2)}$. 
The functions
\beq
 \varphi^{+}_{n_1,n_2}(z) = \sqrt{dz}     \; \; 
\varphi_{n_1,n_2}(z) & \in & C^\infty \left( K_0^{1/2} \otimes W_{(\nu_1,\nu_2)} \right) 
\eeq
form a complete system of eigenfunctions for the operator $(-4 \partial\bar\partial)_{(\nu_1,\nu_2)}$
with the eigenvalues
\[
 \left(\frac{\pi}{ \im\tau \, |\, \pmb{\omega} |}\right)^2 \; \left\vert \left( n_1 + \frac{1-\nu_1}{2}\right) \tau - \left( n_2 + \frac{1-\nu_2}{2}\right) 
 \right\vert^2 \;.
\]
Using the K\"ahler form one can identify $\sqrt{dz}\otimes d\bar{z}$ with $\sqrt{d\bar{z}}$. 
The $\zeta$-function
\begin{equation}
\label{zeta}
 \zeta_{(\nu_1,\nu_2)}(s) = \sum_{n_1,n_2} \frac{1}{\left\lbrack \left( n_1 + \frac{1-\nu_1}{2}\right)^2 \im^2\tau + \left( \left( n_1 + \frac{1-\nu_1}{2}\right)
\re\tau - \left( n_2 + \frac{1-\nu_2}{2}\right)\right)^2 \right\rbrack^s} 
\end{equation}
is absolutely convergent for $\re s>1$. When $(\nu_1,\nu_2)=(1,1)$ it is understood that the summation does not 
include $n_1=n_2=0$.  The function $\zeta_{(\nu_1,\nu_2)}$ is well-defined and has
a meromorphic extension to $\mathbb{C}$ and 0 is not a pole.
The regularized determinant of $-4\,\partial\bar\partial$ is defined by setting
\begin{eqnarray*}
 \ln \det (-4 \partial\bar\partial)_{(\nu_1,\nu_2)} & := & 
 - \left[ \dfrac{1}{\left(\frac{\pi}{ \im\tau \, |\, \pmb{\omega}|}\right)^{2s}}  \; \zeta_{(\nu_1,\nu_2)}(s) \right]' \;.
\end{eqnarray*}
It follows that
\begin{eqnarray*}
 \ln \det (-4 \partial\bar\partial)_{(\nu_1,\nu_2)}
& = & - \zeta'_{(\nu_1,\nu_2)}(0) + \ln\left( \frac{\pi}{ \im\tau \, |\,\pmb{\omega}|} \right)^2 \; \zeta_{(\nu_1,\nu_2)}(0) \;.
\end{eqnarray*}
It was shown in \cite{RS} that $\zeta(0)=0$ for $(\nu_1,\nu_2)\not =(1,1)$, 
and $\zeta(0)=-1$ for $(\nu_1,\nu_2) =(1,1)$. It follows that
\begin{subequations}
\label{determinant_ker}
\begin{align}
 \sideset{}{'}\det (-4 \partial \bar\partial)_{(1,1)} & =  \frac{4 \,\im^2(\tau)  \, |\, \pmb{\omega}|^2}{(2\pi)^2} \; |\pmb{\eta}(\tau)|^4 \;,\\
\label{determinant_gen}
 \det (-4\partial \bar\partial)_{(\nu_1,\nu_2)} & =  \left| \frac{\vartheta_{\nu_1 \nu_2}(\tau)}{\pmb{\eta}(\tau)} \right|^2
= 
e^{-\frac{2\pi}{\im\tau} \, \left(\im \check{\xi} \right)^2} \; \left| 
\frac{\vartheta\left( \check{\xi} \big\vert \tau \right)}{\pmb{\eta}(\tau)} \right|^2
\;,
\end{align}
\end{subequations}
where the Dedekind $\pmb{\eta}$-function and the Jacobi $\vartheta$-function $\vartheta(v|\tau)=\vartheta_{00}(v|\tau)$  are 
given by
\begin{subequations}
\begin{align}
\pmb{\eta}(\tau) & = e^{\frac{\pi i \, \tau}{12}} \; \prod_{n=1}^\infty \left(1- e^{2\pi i n \, \tau}\right) \;, \\
 \vartheta_{ab}(v|\tau)  &= \sum_{n \in \mathbb{Z}} \; \exp\left[ i\pi \left(n+\frac{a}{2}\right)^2 \tau 
+ 2\pi i \left(n+\frac{a}{2}\right)\left(v+\frac{b}{2}\right) \right]
 \;.
\end{align}
\end{subequations}

\section{The topological one-loop string amplitude}
\label{string_amplitude}
For each smooth fiber $E_u$ of the fibration $\Z \to \UP$ we have $\dim H^1(E_u)=2$. Since a base point is given in each fiber by the section $S$, we can 
choose a symplectic basis $\lbrace \alpha_u, \beta_u \rbrace$ of the homology $H_1(E_u)$ with respect to the intersection form, called a {\it homological marking}
consisting of the A-cycle and B-cycle.
We cannot define $\alpha_u, \beta_u$ globally over $\UP$. The cycles are transformed by monodromies around the points
with singular fibers. However, we can define globally an analytical marking.
An {\it analytical marking} is a choice  of a non-zero one-form on each smooth fiber $E_u$. 
We choose the analytical marking that identifies the canonical differential $dx/y$
(where $(x,y)$ are the Weierstrass coordinates in Eq.~(\ref{jacobian}) on the fiber $E_u$)
with the holomorphic one-form $dz$.
Given the elliptic surface $\Z \to \mathrm{UP}$
and the analytic marking we associate to it a holomorphic symplectic two-form \cite{Freed}, given by
\begin{equation}
\label{holomorphic_two_form}
 \pmb{\lambda} = du \wedge \frac{dx}{y} \;.
\end{equation}
Using the two-form $\pmb{\lambda}$ the period integrals of the elliptic fiber $E_u$ over $u$ with periods $2\, \pmb{\omega}, 2\, \pmb{\omega'}$ can be written as follows
\beq
 \int_{\alpha_{u}} \pmb{\lambda} = 2\, \pmb{\omega} \, du \;, \qquad \int_{\beta_{u}} \pmb{\lambda} =  2\, \pmb{\omega'} \, du  \;.
\eeq
Then, there is a globally defined, real closed non-vanishing two-form form $\pmb{\Omega}$ on $\UPo$
\begin{equation}
\label{kaehler_metric}
 \pmb{\Omega} = \int_{E_u} \pmb{\lambda} \wedge \pmb{\overline{\lambda}} = 8 \, i \, \im \tau \; \vert \, \pmb{\omega} \vert^2 \; du\wedge d\bar{u}\;.
\end{equation}
The $\bar\partial$-operator along the fiber $E_u=\pi^{-1}(u)$ is the operator
\begin{equation}
\label{dbar}
 \bar{\partial}: \Omega^{0,0}(E_u) \to \Omega^{0,1}(E_u) \;.
\end{equation}
Its adjoint will be denoted by $\bar{\partial}^*$. We have the following lemma:
\begin{lemma}
\label{regularized_determinant}
The regularized determinant $\sideset{}{^\prime}\det (\bar\partial^* \bar\partial)$ of the Laplace operator along the fiber
is a smooth function on $\UP$ given by
\begin{equation}
\label{det_dbar}
 \sideset{}{^\prime}\det (\bar\partial^* \bar\partial) = \sideset{}{^\prime}\det 
\Delta= \frac{\operatorname{vol}(E_u)^2}{(2\pi)^4} \; \left| \pmb{\Delta}^{\frac{1}{12}}
\right|^2 
\end{equation}
where $\pmb{\Delta}$ is the modular discriminant of the elliptic fiber $E_u$.
\end{lemma}
\begin{proof}
It follows from Eq.~(\ref{determinant_ker}) that
\begin{equation*}
 \sideset{}{^\prime}\det (\bar\partial^* \bar\partial) = \sideset{}{^\prime}\det \Delta
= \frac{\operatorname{vol}(E_u)^2}{(2\pi)^2} \; \left|\frac{\pmb{\eta}^2(\tau)}{2 \,\pmb{\omega}}\right|^2 \;.
\end{equation*}
The discriminant $\pmb{\Delta}$ of the elliptic curve $E_u$ is given by $\pmb{\Delta} = (2\pi)^{12} \, \frac{\pmb{\eta}^{24}(\tau)}{(2\,\pmb{\omega})^{12}}$.
Finally, we have $\operatorname{vol}(E_u)=4 \,\im\tau \, |\, \pmb{\omega}|^2$.
\end{proof}
\noindent
We write $d\mathbf{a} = \pmb{\omega} \, du$ and $d\mathbf{a}_D = \pmb{\omega'} \, du$ such that
\beq
 \tau = \frac{\pmb{\omega'}}{\pmb{\omega}} = \frac{d\mathbf{a}_D}{d\mathbf{a}} \;.
\eeq
This notation should not suggest that $d\mathbf{a}$ is integrable, i.e.,
that there is a globally defined function $\mathbf{a}$. 
On every open set $U \subset \UP$, we can 
integrate and find holomorphic functions $(\mathbf{a},\mathbf{a}_D)$ such that on $U$ we have
$\pmb{\omega} = \frac{d\mathbf{a}}{du}$ and $\pmb{\omega'} = \frac{d\mathbf{a}_D}{du}$.
The K\"ahler metric (\ref{kaehler_metric}) becomes $\pmb{\Omega}= 8 i \, \im\tau \, d\mathbf{a}\wedge d\overline{\mathbf{a}}$.
The following lemma was proved in \cite{Freed}:
\begin{lemma}\label{LCconnection}
The Levi-Civita connection $\nabla^{\operatorname{LC}}$ on $\UPo$ is given by
\begin{equation}
\label{LC_connection}
 \nabla^{\operatorname{LC}} \frac{\partial}{\partial \mathbf{a}} = - \frac{i}{2 \, \im\tau} \, d\tau  \otimes \frac{\partial}{\partial \mathbf{a}} \;, \qquad
 \nabla^{\operatorname{LC}} \frac{\partial}{\partial \mathbf{\bar{a}}} = \frac{i}{2 \, \im\tau} \, d\bar{\tau}  \otimes \frac{\partial}{\partial \mathbf{\bar{a}}} \;.
\end{equation}
The scalar curvature of the Levi-Civita connection is
\begin{equation}
\label{scalar_curvature}
 S  =  \frac{1}{8 \, \im^3 \tau} \; \left\vert \frac{\partial \tau}{\partial \mathbf{a}} \right\vert^2 \;.
\end{equation}
\end{lemma}
\begin{proof}
Let $\pi^{(1,0)} \in \Omega^{(1,0)}(T_{\mathbb{C}}\UP)$ be the projection onto the $(1,0)$ part of the complexified tangent bundle. 
$\pi^{(1,0)}$ is a one-form with values in the tangent bundle $T(\UP)$.
The Levi-Civita connection $\nabla^{\operatorname{LC}}$ on $\UPo$ is defined by
\beq
  d \, \left[ \pmb{\Omega}\left( \frac{\partial}{\partial \mathbf{a}} \, , \, \frac{\partial}{\partial \mathbf{\bar{a}}} \right)\right]
= \pmb{\Omega}\left( \nabla^{\operatorname{LC}} \, \frac{\partial}{\partial \mathbf{a}} \, , \, \frac{\partial}{\partial \mathbf{\bar{a}}} \right)
 - \pmb{\Omega}\left( \frac{\partial}{\partial \mathbf{a}} \, , \, \nabla^{\operatorname{LC}} \, \frac{\partial}{\partial \mathbf{\bar{a}}} \right) \;.
\eeq
It is the unique connection which is compatible with the metric and the complex structure.
It follows that the Levi-Civita connection satisfies $\pmb{\Omega}\left( \pi^{(1,0)}, \nabla^{\operatorname{LC}} \pi^{(1,0)}\right) = 0$.
Eqns.~(\ref{LC_connection}) follow. We find
\beq
 d^{\nabla^{\operatorname{LC}}} \left( \nabla^{\operatorname{LC}} \frac{\partial}{\partial \mathbf{a}} \right) & = & - \frac{i}{2}\, d\left( \frac{1}{\im \tau} \right) \wedge d\tau \otimes 
\frac{\partial}{\partial \mathbf{a}} - \frac{i}{2 \, \im\tau} \, d\tau  \wedge  \nabla^{\operatorname{LC}} \frac{\partial}{\partial \mathbf{a}} \\
& = & - \frac{1}{4 \, \im^2\tau} \, d\tau\wedge d\bar\tau \otimes 
\frac{\partial}{\partial \mathbf{a}} \;.
\eeq
The Riemannian curvature $R$ of the K\"ahler metric $\pmb{\Omega}$ is
\beq
 R = \left( \begin{array}{cc} R^{\;\; a}_{a \;\; a\bar{a}} & 0 \\ 
 0 & R^{\;\; \bar{a}}_{\bar{a} \;\; a\bar{a}}  \end{array} \right)
\; d\mathbf{a}\wedge d\overline{\mathbf{a}} 
\eeq
where
\beq
 R^{\;\; a}_{a \;\; a\bar{a}} = d\mathbf{a} \left\lbrace \; d^{\nabla^{\operatorname{LC}}}  \left( \nabla^{\operatorname{LC}} \frac{\partial}{\partial \mathbf{a}} \right)(\partial_{\mathbf{a}} \,,\, 
\partial_{\mathbf{a}} ) \; \right\rbrace = - \frac{1}{4 \, \im^2\tau} \left\vert \frac{\partial\tau}{\partial \mathbf{a}} \right\vert^2\;.
\eeq
Pulling down the summation index using the metric $\pmb{\Omega} = i \,h_{a\bar{a}} \, d\mathbf{a}\wedge d\overline{\mathbf{a}}$ we find
\beq
 R_{a \bar{a} a\bar{a}} & = &  h_{a\bar{a}} \; R^{\;\; a}_{a \;\; a\bar{a}} 
= - \frac{2}{\im \tau} \; \left\vert \frac{\partial \tau}{\partial \mathbf{a}} \right\vert^2 \;.
\eeq
The scalar curvature is obtained by contracting the summation indices using the inverse K\"ahler metric $\Omega$. We obtain
\beq
 S & = & -4 \, \left( h^{a\bar{a}} \right)^2  \; R_{a \bar{a} a\bar{a}} = \frac{1}{8 \, \im^3 \tau} \; \left\vert \frac{\partial \tau}{\partial \mathbf{a}} \right\vert^2 \;.
\eeq
\end{proof}
\noindent
We obtain the following proposition:
\begin{proposition}
The function $F^{(1)} = - \frac{1}{2} \ln \sideset{}{^\prime}\det \Delta$ is a smooth function on $\UPo$
and satisfies the equation
\begin{equation}
\label{anomaly}
 \Delta_{\UP} \, F^{(1)} = S \;,
\end{equation}
where $\Delta$ is the Laplace operator along the fiber of $\Z \to \UP$, $S$ is the scalar curvature of the 
K\"ahler metric $\pmb{\Omega}$ on $\UPo$, and $\Delta_{\UP}$ is the Laplace-Beltrami operator 
\begin{equation}
  \Delta_{\UP} = \frac{1}{\im\tau} \,  \partial_{a} \partial_{\bar{a}}  \;.
\end{equation}
\end{proposition}
\begin{proof}
The proof follows from Lemma \ref{regularized_determinant} and Lemma \ref{LCconnection}.
\end{proof}
\begin{remark}
Eq.~(\ref{anomaly}) is the anomaly equation of the one-loop topological string amplitude derived in \cite{BCOV}, i.e., 
\begin{equation}
  \partial_{a} \partial_{\bar{a}}   \, F^{(1)} = \frac{1}{8  \, \im^2\tau}  \, \left\vert \frac{\partial \tau}{\partial \mathbf{a}} \right\vert^2  \;.
\end{equation}
\end{remark}

\section{The vertical $\bar{\partial}$-operator on $\Z \to \UP$}
\label{relative_determinant0}
\label{global_anomaly}

In Section~\ref{single_curve}, we computed the regularized determinant of the $\bar{\partial}$-operator
on an elliptic curve. When the elliptic curve $E_u$ is varied in an elliptic surface, we obtain the determinant 
line bundle $\DET \bar{\partial} \to \UP$ of the $\bar{\partial}$-operator along the fiber of $\Z \to \UP$. 
The determinant line bundle is the holomorphic line bundle
\begin{equation}
\label{determinant_cohomology}
 \DET \bar{\partial} \to \UPo \, \quad \text{with fibers} \quad \left(\DET \bar{\partial}\right)_{u} 
 =   H^{0,0}(E_u, \mathbb{C})^{-1} \otimes  H^{0,1}(E_u, \mathbb{C})  \;.
\end{equation}
There is a factorization of the determinant line bundle 
$\mathcal{L}=\DET \overline{\partial}$ as the tensor product of $\mathcal{L}'$  and $\mathcal{H}$, 
corresponding to the non-zero and zero eigenvalues respectively. 
The line bundle $\DET \bar\partial$ can be identified with $\mathcal{H}$ since
\beq
 H^{0,0}(E_u, \mathbb{C})   \cong   \ker \bar\partial \;,\quad
 H^{0,1}(E_u, \mathbb{C})   \cong   \ker \bar\partial^*  \;.
\eeq 
The bundle $\mathcal{L}'$ has a holomorphic section $\sideset{}{'}\det 2\bar\partial$ that
determines an isomorphism $\DET \bar{\partial} \cong \mathcal{H}$ \cite{Atiyah2}.
The isomorphism does not preserve the metric or connection. Bismut and Freed defined the smooth metric 
\begin{equation}
\label{determinant_metric}
 \Vert.\Vert_Q :=  (2\pi)^2 \, \sqrt{\sideset{}{^\prime}\det \Delta} \; \; \Vert.\Vert_{L^2} 
\end{equation}
on $\DET \bar\partial$ and determined its unitary connection. Since $\omega_{Z/\UP}$ is equipped with a Hermitian $C^\infty$-metric,
the Quillen metric is a Hermitian $C^\infty$-metric on the holomorphic fibers $(\DET \bar{\partial})_{u}$.
The curvature of this connection is called the local anomaly in physics.
It follows:
\begin{lemma}
\label{laplace}
The determinant line bundle $\DET \bar{\partial} \to \UPo$ 
is flat. $\sigma=(dz)^{-1}$ is a non-vanishing holomorphic section of $\DET \bar{\partial}$ 
with
\beq
 \Vert \sigma \Vert_Q = |\pmb{\Delta} |^\frac{1}{12} \;.
\eeq
$\sigma^* = dz$ is a non-vanishing holomorphic section of the dual bundle $(\DET \bar{\partial})^*\to \UPo$
with $\Vert \sigma^* \Vert_{Q^*} = |\pmb{\Delta} |^{-\frac{1}{12}}$.
\end{lemma}
\begin{proof}
The flatness follows from the curvature formula of Bismut and Freed 
(\cite{BismutFreed} or \cite[Prop.~5.14]{Atiyah2}):
\beq
 c_1\left( \DET \bar{\partial} \right) = - \int_{E_u} \frac{c_1^2\left(\omega_{\Z/\UP}\right)}{6} \;.
\eeq
It follows from Remark \ref{canonical_bundle} that $c^2_1( \omega_{\Z/\UP} )=0$.

$H^{0,1}(E_u)$ and  $H^{1,0}(E_u)$ are Serre duals. Thus, we have $\mathcal{H} \cong  [H^{1,0}(E_u) \otimes H^{0,0}(E_u)]^{-1}$.
The kernel consists of the constant function $\phi=1$ with $\Vert \phi 
\Vert^2 = \operatorname{vol}(E_u)$. 
By Serre duality we identify the cokernel $\ker \bar\partial^*$ with the dual of the space of holomorphic one-forms. 
Thus, the cokernel is spanned by the section $(dz)^{-1}$ and $\Vert dz \Vert^2 = \operatorname{vol}(E_u)$. 
Using Eq.~(\ref{determinant_ker}) we obtain for the Quillen norm of the section $(dz)^{-1}$
\beq
 \Big\Vert (dz)^{-1} \Big\Vert^2_Q = (2\pi)^4 \; \frac{\det (-4\partial \bar\partial)_{(1,1)}}{\Vert \phi\Vert^2 \, 
\Vert dz \Vert^2} = 
 \left\vert\pmb{\Delta}\right|^\frac{2}{12} \;.
\eeq
It is possible to factorize the right hand side holomorphically in $\tau$. 
We use the Quillen metric to obtain a smooth section $\sigma^{\#}$ of the dual bundle $(\DET \bar{\partial})^*$ 
by setting
\beq
 \sigma^{\#} = g_Q( \sigma, \bullet) = |\pmb{\Delta}|^\frac{1}{6} \, dz \;,
\eeq
and $\Vert \sigma^{\#} \Vert_{Q^*} = |\pmb{\Delta}|^{1/12}$ by definition. 
The claim follows.
\end{proof}
\begin{remark}
\label{change_of_variables}
Under the change of the coordinate chart from $[u:1]$ to $[1:v]$ on $\CP1$ by $u=-1/v$
the holomorphic differential transforms as $dz_u = \frac{dx_u}{y_u} = - \, v \, \frac{dx_v}{y_v} = - \, v \, dz_v$.
The Quillen metric is compatible with a change of coordinates since
\beq
 \Big\Vert (dz_v)^{-1} \Big\Vert_Q^2 = |v|^2 \; \Big\Vert (dz_u)^{-1} \Big\Vert_Q^2 
 = |v|^2 \, \Big|\pmb{\Delta}^{\frac{1}{12}}_{u}\Big|^2 =  \Big|\pmb{\Delta}^{\frac{1}{12}}_v\Big|^2 \;.
\eeq
\end{remark}
\begin{lemma}
\label{bismut_freed}
The Bismut-Freed connection on $\DET \bar\partial \to \UPo$ is flat and given by
\beq
 \nabla^{\operatorname{BF}} \sigma &= & \frac{1}{12} \; \frac{\partial \pmb{\Delta}}{ \pmb{\Delta}} \otimes \sigma \;,
\eeq
and $\nabla^{\operatorname{BF} \, (0,1)} = \bar{\partial}$. 
\end{lemma}
\begin{proof}
It follows for the Quillen metric in the coordinate chart $[u:1] \in \CP1$
\beq
\frac{d \Vert \sigma \Vert^2_Q}{\Vert \sigma \Vert^2_Q} = \frac{1}{12}\, \frac{ \partial \ln \pmb{\Delta}}{\partial u} \; du 
+ \frac{1}{12}\, \frac{ \partial \ln \bar{\pmb{\Delta}}}{\partial \bar{u}} \; d\bar{u}\;.
\eeq
\end{proof}
Although $\DET \bar{\partial}$ is flat and hence the local anomaly vanishes, there is a global anomaly arising from monodromy 
around the non-contractible closed loops. These are the non-contractible loops $(\gamma_n)_{n=1}^{N_f+2}$ encircling the 
nodes at $(u^*_n)_{n=1}^{N_f+2}$ clockwise, and $\gamma_{\infty}$ encircling the cusp at infinity counterclockwise. 
\begin{lemma}
\label{asymptotics}
There exist constants in $\mathbb{R}^+$ such that
\beq
  \Vert \sigma \Vert_Q & \sim & c_{n} \,\; |u - u^*_n|^{\frac{1}{12}} \qquad ( u \to u^*_n) \;, \\
  \Vert \sigma \Vert_Q & \sim & c_{\infty} \, |v|^{\frac{10-N_f}{12}} \quad\qquad ( v \to 0) \;.
\eeq
\end{lemma}
\begin{proof}
The proof follows from Lemma~\ref{laplace}, Remark~\ref{change_of_variables}, and the fact that $\pmb{\Delta}$
has a simple zero at each node $u_n$ and is a polynomial of degree $N_f+2$ in $u$.
\end{proof}
We denote the holonomy of the determinant section of 
$\DET \bar{\partial}$ on the boundary circle $\gamma_n$ around $u_n$ by $\exp( -\frac{i \pi}{2} \, \eta_{\bar{\partial}}[\gamma_n])$. 
The following lemma is an immediate consequence of Lemma \ref{asymptotics}:
\begin{lemma}
\label{monodromy_and_eta}
The holonomy of the section $\sigma$ on a cycle $\gamma_i$ is $\exp\left(-\frac{\pi i}{2} \, \eta_{\bar{\partial}}[\gamma_i]\right)$
with
\beq
 \eta_{\bar{\partial}}[\gamma_{n}]     \equiv  - \frac{1}{3} \mod{4} \;, \quad 
 \eta_{\bar{\partial}}[\gamma_{\infty}] \equiv  -  \frac{10-N_f}{3} \mod{4}\;.
\eeq
\end{lemma}
Since we have restricted ourselves to the case where the fibration $\Z \to \UP$ has no surface singularities,
the singular fibers $E_u = \pi^{-1}(u)$ at $u=u^*_1,\dots, u^*_{N_f+2}$ are nodal curves. A current $\delta(u - u^*_n)$ on 
$\Z\to \UP$ is defined by saying that for every differential form $\alpha$ on $\Z$ with compact support, the equality
\beq
  \int_Z \delta(u-u^*_n) \, du \wedge d\bar{u}\wedge \alpha = \int_{E_{u^*_n}} \alpha  
\eeq  
holds. The first Chern class $c_1(  \DET \bar\partial , \Vert .\Vert _Q)$ is defined as a current
of type $(1,1)$ by
\begin{eqnarray}
\label{c1_index}
  c_1\left(  \DET \bar\partial, \Vert. \Vert_Q\right) := \frac{1}{2\pi i} \, \partial \bar{\partial} \log \Vert\sigma\Vert_Q^2 
 \end{eqnarray}
where $\sigma$ is a local non-vanishing holomorphic section.
\begin{proposition}
\label{c2_BPSline}
The determinant line bundle with the Quillen metric $\Vert.\Vert_Q$ of the $\bar{\partial}$-operator along the fiber of $\, \Z \to \UP$ extends to a 
holomorphic line bundle $\DET \bar{\partial} \to \UP$  
with curvature
\begin{eqnarray}
\label{curvature}
 c_1\left( \DET \bar{\partial} \right)  & =  & - \sum_{n=1}^{N_f+2} \frac{1}{12} \, \delta(u-u^*_n) \, du \wedge d\bar{u}\;.
\end{eqnarray}
\end{proposition}
\begin{proof}
The definition of the holomorphic determinant line bundle in Eq.~(\ref{determinant_cohomology})
can be extended across the singular fibers of an elliptic fibration using the results of
Knudsen and Mumford \cite{BismutB1}. The Quillen metric is smooth by Lemma \ref{det_dbar}.
Let $f$ be a function that is differentiable in the disc $D_{\epsilon}$ with $|u|<\epsilon$. Suppose further that $f$ and its derivatives
with respect to $\bar{u}$ are bounded on the disc. Let $\mathcal{T}$ denote the function
\begin{eqnarray}
\label{curvature1}
 \Big(\mathcal{T} f\Big)(u,\bar{u}) = \frac{1}{2\pi i} \iint_{|w|< \epsilon} f(w,\bar{w}) \; \; \frac{dw \, d{\bar{w}}}{u - w} \;.
\end{eqnarray}
It is well-known \cite{Nickerson} that the linear operator $\mathcal{T}$ is differentiable and admissible on $D_{\epsilon}$ and satisfies $\bar{\partial}(\mathcal{T} f) = f$.
In this sense, we write
\beq
 \frac{1}{2\pi i} \, \frac{\partial}{\partial \bar{u}} \, \frac{1}{u-w} =  \delta(u-w) \;.
\eeq
Since $\bar\partial \partial + \partial \bar\partial =0$,
Eq.~(\ref{curvature}) follows from Eq.~(\ref{c1_index}) and the application of Lemma~\ref{asymptotics}.
\end{proof}
\begin{remark}
Corollary~\ref{c2_BPSline} shows that we can include the nodal fibers of $\Z \to \UP$ 
when considering the determinant line bundle of the $\bar{\partial}$-operator along the fiber. The contributions $\eta_{\bar{\partial}}[\gamma_{n}]$
to the global anomaly of the determinant line bundle around the nodal fibers
 are then viewed as current contributions of type $(1,1)$ to the
first Chern class of the extended determinant line bundle.
\end{remark}
The results in \cite[Thm.~2.1]{BismutB1} and \cite{BismutB2} describe the more general situation
of the $\bar{\partial}$-operator coupled to a holomorphic vector bundle $V \to \Z$.  There, the authors prove that
\begin{equation}
 \label{index2}
\begin{split}
c_1\left(  \DET \bar\partial_{V}, \Vert. \Vert_Q\right)
= & - \int_\pi \Big[ \operatorname{ch}(V) \wedge 
 \operatorname{Todd}(\omega_{Z/\UP} , \Vert .\Vert ) \Big]_{(4)} \\
& -  \sum_{n} \frac{\operatorname{rank}(V)}{12} \,
\delta(u-u^*_n) \, du \wedge d\bar{u}\;.
\end{split}
\end{equation}
For $\operatorname{Todd}(\omega_{Z/\UP} , \Vert .\Vert ) = 1 + c_1(\omega_{Z/\UP})/2$ (since $c^2_1(\omega_{Z/\UP})=0$),
and $\operatorname{ch}(V)= 1$ we obtain Eq.~(\ref{curvature}).

\section{Extending $\bar{\partial}$ to nodal curves}
\label{singular_dbar}

Let us restrict the fibration $\Z \to \UP$ to a small neighborhood $U$ of a point $u^*$ whose singular fiber is a nodal curve. 
We identify any smooth elliptic fiber $E_u$ of the fibration $\Z \to U$ with the complex plane (with complex coordinate $z$) 
modulo the action of the lattice generated by the periods $2 \, \pmb{\omega}, 2 \, \pmb{\omega'}$. For $\tau = \frac{\pmb{\omega'}}{\pmb{\omega}}$
we set $q=\exp{2\pi i \tau}$. After a suitable $\mathrm{SL}(2,\mathbb{Z})$ transformation we can assume
that as we approach the node for $u \to u^*$ we have $\im\tau \to \infty$, 
$q \to 0$. By making the neighborhood $U$ smaller if necessary, we assume that $|q|<1$ uniformly in $U$.

Next, we consider the annulus $\pmb{\mathrm{ann}}(r_1,r_2)$ in $\mathbb{C}$ (with complex coordinate $W$) with inner radius $0<r_1<1$ and outer radius $r_2=1/r_1$.
We also set $Z=1/W$ for $W\not=0$. If we set $r_1= |q|^{1/2}$, then the annulus is covered by the two charts $r_1<|W|\le 1$ and $r_1 < |Z| \le 1$.
The inner and outer radius of the annulus are identified to obtain a torus using $Z \, W = q$. This is the local description of a 
compact Riemann surface near a node used in \cite{SeeleyS} when applied to a torus. For $q\to 0$ we get a singular surface envisioned as sphere with two 
points $\lbrace Z=0\rbrace$ and $\lbrace W=0\rbrace$ identified.

\noindent
The map $W=r_2 \; \exp{\left( 2\pi i\frac{z}{2\,\pmb{\omega}} \right)}$ identifies the fundamental domain for the torus  $E_u \cong \mathbb{C}/\langle 2\, \pmb{\omega}, 2\, \pmb{\omega'}\rangle$ 
with its vertical edges parallel to $2 \,\pmb{\omega'}$ identified with the annulus $\pmb{\mathrm{ann}}(r_1,r_2)$. 
With respect to the metric $g=dz.d\bar{z}$ we have $\operatorname{vol}(E_u)=4 \,\im\tau \, |\, \pmb{\omega}|^2$
and 
\begin{equation}
\label{metric_g}
 g = dz . \, d\bar{z} = \left( \dfrac{|\pmb{\omega}|}{\pi |W|} \right)^2 \; dW.\,d\overline{W} \;.
\end{equation}
The following lemma computes the regularized determinant for the Laplace-Beltrami operator on the annulus 
with respect to the metric (\ref{metric_g}) and with Dirichlet boundary conditions, i.e., for eigenfunctions vanishing 
on the outer and inner radius:
\begin{lemma}
The regularized determinant for the Laplace operator $\Delta=-4\partial_z \bar{\partial}_z$ on the annulus 
$\pmb{\mathrm{ann}}(r_1,r_2)$ with Dirichlet boundary conditions is given by
\begin{equation}
\label{det_Dirichlet}
\sideset{}{_D}\det \Delta 
= \frac{\im\tau}{2\pi}  \; \left|\,\pmb{\eta}(\tau)\right|^2  \;.
\end{equation}
\end{lemma}
\begin{proof}
The eigenfunctions  $\psi_{n_1,n_2}(\xi^1,\xi^2)$ of the Laplace-Beltrami operator with Dirichlet boundary conditions are similar to the 
eigenfunctions in section \ref{single_curve}. In the case $(\nu_1,\nu_2)=(1,1)$ we have
\begin{equation*}
 \varphi_{n_1,n_2}(\xi^1,\xi^2) 
 =   \exp\left[2\pi i n_1 \left( \xi^1 - \frac{\re\tau}{\im\tau} \xi^2\right) \right]
\;  \left[ \cos\left(\frac{2\pi n_2}{\im\tau} \xi^2\right) + i \, \sin\left(\frac{2\pi n_2}{\im\tau} \xi^2\right) \right]\;.
\end{equation*}
The Dirichlet boundary conditions are $\psi_{n_1,n_2}(\xi^1,0) = \psi_{n_1,n_2}(\xi^1,\im\tau) =0$. Hence, we have $n_2>0$
and
\begin{equation*}
 \psi_{n_1,n_2}(\xi^1,\xi^2) = \exp\left[2\pi i n_1 \left( \xi^1 - \frac{\re\tau}{\im\tau} \xi^2\right) \right]
\;  \sin\left(\frac{2\pi n_2}{\im\tau} \xi^2\right) \;.
\end{equation*}
Therefore, the zeta-function $\zeta_D(s)$ for the Laplace operator with Dirichlet boundary conditions is
\begin{equation*}
\begin{split}
 \zeta_D(s) = \sum_{n_1,n_2>0} \frac{1}{\left\lbrack n_1^2 \, \im^2\tau + \left( n_1 \re\tau -  n_2 \right)^2 \right\rbrack^s} = \frac{1}{2} \left( \zeta_{(1,1)}(s) - \dfrac{2}{|\tau|^{2s}} \, \zeta(2s) \right)\;,
\end{split}
\end{equation*}
where $\zeta_{(1,1)}(s)$ was defined in Eq.~(\ref{zeta}) and $\zeta(2s)=\sum_{n=1}^\infty n^{-2s}$ is the Riemann zeta function. 
Thus, we have
\begin{equation*}
\begin{split}
 \zeta_D(0)  = \frac{1}{2} \, \left( \zeta_{(1,1)}(0) - 2 \, \zeta(0) \right) \;,\quad
  \zeta'_D(0)  = \frac{1}{2} \, \left( \zeta'_{(1,1)}(0) + 4 \, \ln{|\tau|} \, \zeta(0) - 4\, \zeta'(0) \right) \;.
\end{split}
\end{equation*}
Eq.~(\ref{det_Dirichlet}) then follows from $\zeta(0)=-1/2$, $\zeta'(0)=-\ln{(2\pi)}/2$, and Eq.~(\ref{determinant_ker}).
\end{proof}
Instead of the Laplacian $\Delta = - 4 \partial_z\bar{\partial}_z$ Seeley and Singer \cite{SeeleyS} use the Laplace-Beltrami operator for the flat metric on the annulus.
The flat metric $\widehat{g}$ on the annulus is obtained from the metric $g$ in Eq.~(\ref{metric_g}) by a change in the conformal gauge, i.e.,
\begin{equation}
\label{metric_ghat}
 \widehat{g} = dW.\,d\overline{W} = e^{2\, \Phi} \; g
\end{equation}
with $\Phi = \ln{(\pi\, |W|/|\pmb{\omega}|)}$. It follows that the Laplace-Beltrami operator for the flat metric is given by
\begin{equation}
\label{Laplace_ghat}
 \widehat{\Delta} = - 4 \, \partial_W \, \bar{\partial}_W =  e^{-2\,\Phi}  \, \Delta\;.
\end{equation}
We use the results of \cite{Weisberger} to derive the following lemma:
\begin{lemma}
The regularized determinant of the Laplace operator $\widehat{\Delta}$ on the annulus $\pmb{\mathrm{ann}}(r_1,r_2)$ with Dirichlet boundary conditions is
given by
\begin{equation}
\label{laplace_ghat}
\sideset{}{_D}\det \widehat{\Delta}
= \frac{\im\tau}{2\pi}  \; \left|\,\pmb{\eta}(\tau)\right|^2 \; |q|^{\frac{1}{6}}\;.
\end{equation}
\end{lemma}
\begin{proof}
The Gaussian curvature as well as the geodesic curvature on the boundary vanish for $g$. 
Similarly, the Gaussian curvature and the average geodesic curvature on the boundary vanish for $\widehat{g}$. Then, \cite[Eq.~(3)]{Weisberger} implies
\begin{equation}
 \sideset{}{_D}\det \Delta = \sideset{}{_D}\det \widehat{\Delta} \; \; \exp\left(\frac{L}{6\pi}\right) 
\end{equation} 
where $L= \frac{1}{2} \int_{\mathrm{ann}(r_1,r_2)} \mathrm{vol}_{g}  \; g^{ab} \, (\partial_a\Phi)\, (\partial_b\Phi)$.
A calculation shows that
\begin{equation*}
 L= \frac{i}{4} \int_{\mathrm{ann}(r_1,r_2)} \frac{dW\wedge d\overline{W}}{|W|^2} = \pi \ln\left(\frac{r_2}{r_1}\right) = 2 \pi^2 \im\tau 
\end{equation*}
and $ \exp{L/6\pi} = |q|^{(-1/6)}$.
We obtain
\begin{equation}
 \sideset{}{_D}\det \Delta  = \sideset{}{_D}\det \widehat{\Delta} \; \, |q|^{-\frac{1}{6}} \;.
\end{equation} 
\end{proof}

\noindent
The application of the result of Seeley and Singer \cite{SeeleyS}  yields the following proposition:
\begin{proposition}
In a small neighborhood $U\subset \UP$ of a point $u^*$ with nodal fiber $E_{u^*}$ such that
$q=\exp(2\pi i \tau)\to 0$ as $u \to u^*$ the family of operators $\lbrace \bar{\partial}_{W,u} \rbrace_{u \in U}$ is a continuous family
and the operator $\bar{\partial}_W$ is well-defined on the singular fiber $E_{u^*}$.
\end{proposition}
\begin{remark}
The limiting Laplace operator of \cite{SeeleyS} is the Laplace operator $- 4 \, \partial_W \, \bar{\partial}_W$ on $\mathbb{C}$ 
(with complex coordinate $W$). The eigenfunction for an eigenvalue $\lambda^2$ with $\lambda>0$ satisfying Dirichlet boundary conditions is 
$\mathrm{J}_n(\lambda \,r) \, \exp(in\theta)$ with $W=r \, \exp(i\theta)$ and $n\in\mathbb{N}$. $\mathrm{J}_n(\lambda \,r)$ is the
Bessel function of the first kind that is regular at $r=0$ and decays as $1/\sqrt{r}$ for $r\to \infty$.
\end{remark}

\section{The vertical signature operator on $\Z \to \UP$}
\label{relative_signature}
Using its complex structure the signature operator on each fiber $E_u$ can be identified with the operator 
\beq
 \mathrm{D}=\overline{\partial} + \overline{\partial}_1: \quad \Omega^{0,0}(E_u) \oplus \Omega^{1,0}(E_u)
 \to \Omega^{0,1}(E_u) \oplus \Omega^{1,1}(E_u) \;.
\eeq
Again, there is a factorization of the determinant line bundle 
$\mathcal{L}=\DET \, \mathrm{D}$ as the tensor product of $\mathcal{L}'$  and $\mathcal{H}$, 
corresponding to the non-zero and zero eigenvalues respectively \cite{Atiyah2}. The bundle $\mathcal{L}'$ has the 
holomorphic section $\sideset{}{'}\det \mathrm{D}$. The fiber of the line bundle $\mathcal{H}$ is
\beq
 \mathcal{H} \cong \Big[ H^{0,0}(E_u) \otimes H^{1,0}(E_u) \Big]^{-1} \otimes \Big[ H^{0,1}(E_u) \otimes H^{1,1}(E_u) \Big]\;.
\eeq
The bundles $H^{0,0}(E_u)$ and $H^{1,1}(E_u)$ can be identified by duality. Similarly, $H^{0,1}(E_u)$
and  $H^{1,0}(E_u)$ are Serre duals. On each fiber $E_u$ multiplication by $dz$ converts $\overline{\partial}$ into 
$\overline{\partial}_1$. Thus, we have $\mathcal{H} \cong  [H^{1,0}(E_u)]^{-2}$
and the determinant line bundle of the operator $\mathrm{D}$ is isomorphic to the two-fold tensor product 
of $\DET\overline{\partial}$. 

We have the following lemma:
\begin{lemma}
\label{local_signature}
\begin{enumerate}
\item[]
\item $\sigma_{\mathrm{D}}=(dz)^{-2}$ is a non-vanishing holomorphic section of $\DET \; \mathrm{D}\to \UPo$ 
with
\beq
 \Vert \sigma_{\mathrm{D}} \Vert_Q = |\pmb{\Delta} |^\frac{1}{6} \;.
\eeq
$\sigma^*_{\mathrm{D}} = (dz)^2$ is a holomorphic section of the dual bundle $(\DET \; \mathrm{D})^*\to \UPo$
with $\Vert \sigma^*_{\mathrm{D}} \Vert_{Q^*} = |\pmb{\Delta} |^{-\frac{1}{6}}$.
\item The flat Bismut-Freed connection on $\DET \; \mathrm{D}\to \UPo$ is given by
\beq
 \nabla^{\operatorname{BF}} \sigma_{\mathrm{D}} &= & \frac{1}{6} \frac{\partial \pmb{\Delta}}{ \pmb{\Delta}} \otimes \sigma_{\mathrm{D}} \;,
\eeq
and $\nabla^{\operatorname{BF} \, (0,1)} = \bar{\partial}$. 
\item The determinant line bundle with the Quillen metric extends to a holomorphic line bundle
$\DET \; \mathrm{D} \to \UP$. The curvature is a current with
\beq
 c_1\left( \DET \; \mathrm{D} \right)  & =  & - \sum_{n=1}^{N_f+2} \frac{1}{6} \, \delta(u-u^*_n) \, du\wedge d\bar{u}\;.
\eeq
\end{enumerate}
\end{lemma}
\begin{proof}
The proof is the same as for Lemma \ref{bismut_freed}.
\end{proof}
\begin{lemma}
\label{trivialization}
The line bundle $(\DET \; \mathrm{D})^{* \, 6} \to \UPo$ is canonically trivial.
\end{lemma}
\begin{proof}
The section $\underline{\sigma}=\pmb{\Delta}^{1/6} \; (dz)^{2}$ of $(\DET \; \mathrm{D})^*\to \UPo$ satisfies $\Vert \underline{\sigma} \Vert_{Q^*}=1$
and is invariant under the action of $\pi_1(\UPo)$ up to a sixth root of unity.
The trivializing, holomorphic, non-vanishing section of $(\DET \; \mathrm{D})^{*\otimes 6}\to \UPo$
is $\underline{\sigma}^{6}$.
\end{proof}
\noindent
It follows from Remark~\ref{change_of_variables2} that $\mathcal{H} \cong T^*\UP$.
Consequently, we can obtain well-defined {\em logarithmic} monodromies for the bundle $(\DET \; \mathrm{D})^{*} \to \UPo$. 
We denote this distinguished choice for the
monodromy by $\eta^0_{\mathrm{D}}[\gamma]$ as opposed to $\eta_{\mathrm{D}}[\gamma]$ which appeared in Lemma \ref{monodromy_and_eta}
and was only determined modulo $4$. 
\begin{lemma}
\label{log_monodromy}
The logarithmic monodromies of the bundle $(\DET \; \mathrm{D})^{*} \to \UPo$ are
\beq
 \eta^0_{\mathrm{D}}[\gamma_n]   = - \frac{2}{3} \;, \qquad
 \eta^0_{\mathrm{D}}[\gamma_{\infty}] =   - \frac{2(10-N_f)}{3} \;.
\eeq
\end{lemma}
\begin{proof}
Under the isomorphism $\mathcal{H} \cong T^*\UP$ the form $(dz)^{2}$ is identified with $du^{-1}$.
Thus, the bundle $(\DET \; \mathrm{D})^{*} \otimes T^*\UP \to \UPo$ has a standard trivialization on $\UPo$ given by $(dz)^2 \otimes du^{-1}$.
In this trivialization the holomorphic section is $\pmb{\Delta}^{1/6}$. The claim follows.
\end{proof}

\begin{lemma}
\label{signature_coeff}
The signature of the elliptic surface $\Z \to \UP$ is
\beq
 \operatorname{sign}(\Z) =  \sum_{n=1}^{N_f+2}  \eta^0_{\mathrm{D}}[\gamma_n] - \frac{1}{2} \eta^0_{\mathrm{D}}[\gamma_{\infty}] + 2 \;.
\eeq
It follows $\operatorname{sign}(\Z)=-N_f$.
\end{lemma}
\begin{proof}
The signature of the rational elliptic surface $\Zc \to \CP1$ in terms of its Chern classes $c_1, c_2$ is 
\beq
 \operatorname{sign}(\Zc) = \frac{c_1^2-2 c_2}{3} \;.
\eeq
The canonical class is minus the fiber class so that $c_1^2=0$ while $c_2$ is the sum of 
the exceptional fibers of the fibration whence
\beq
 \operatorname{sign}(\Zc) = - \frac{2}{3} \sum_{n=1}^{N_f+2} e(E_{u^*_n})  - \frac{2}{3} e(E_{\infty})\;.
\eeq
The elliptic surface $\Zo \to \UPo$ is obtained by cutting out all singular fibers
$\Zo = \Zc - \cup E_{u^*_n} - E_{\infty}$; the elliptic surface $\Z \to \UP$ is obtained by cutting out 
only the singular fiber at infinity whence $\Z = \Zc - E_{\infty}$. Since the singular fibers at $u=u^*_1, \dots, u^*_{N_f+2}$
are nodes they do not contribute to the signature. Hence, we have
\beq
   \operatorname{sign}(\Zo) = \operatorname{sign}(\Z) = \operatorname{sign}(\Zc) - \operatorname{sign}(E_{\infty}) \;,
\eeq
where $E_\infty$ is the singular fiber of $\Zc$ at infinity.
Thus, we obtain
\beq
 \operatorname{sign}(\Z) =  - \frac{2}{3} \sum_{n=1}^{N_f+2} e(E_{u^*_n})  - \frac{2}{3} e(E_{\infty}) - \operatorname{sign}(E_{\infty}) \;.
\eeq
The Euler number $e(E_{u^*_n})$ is equal to the degree of the zero the discriminant assumes at $u=u^*_n$. Therefore,
it follows $-\frac{2}{3} e(E_{u^*_n})=  \eta^0_{D}[\gamma_n]$ and $-\frac{2}{3} e(E_{\infty})=  \eta^0_{D}[\gamma_{\infty}]$.
By Kodaira's classification result it follows that singularities which are not of type 
$I_k$ satisfy $\operatorname{sign}(E_{\infty})= 2 - e(E_{\infty})$. Lemma~\ref{log_monodromy} yields $\operatorname{sign}(\Z)=-N_f$.
\end{proof}

\subsection{Regular singularities and the Riemann-Hilbert problem}
\label{riemann_hilbert}

The definition of the holomorphic determinant line bundle in Eq.~(\ref{determinant_cohomology})
can be extended across the singular fibers of an elliptic fibration using the results of Knudsen and Mumford \cite{BismutB1} to include the higher-rank singularities
of Kodaira type $I_k$ and $I_k^*$. In this section, we allow any Jacobian rational elliptic surface $\Z$ with singular fibers of Kodaira type 
$I_{k_n}$ over $[u^*_n:1]\in \CP1$ (with $1 \le n \le K$ such that 
$\sum k_n = N_f +2$) and a singular fiber of Kodaira type $I^*_{4-N_f}$ over $u=\infty$. 
The following solution to the Riemann-Hilbert problem on $\CP1$ was given by R\"ohrl in terms of differential 
equations with regular singular points \cite{Majima}:

\begin{fact}
\label{roehrl}
\begin{enumerate}\item[]
\item The functor mapping all conjugate classes of one-dimensional representations of $\pi_1(\UPo)$ to the set of isomorphism classes
of flat line bundles over $\UPo$ is an equivalence of categories. \hoehe
\item For a flat line bundle $\mathfrak{E} \to \UPo$ together with the natural connection $d$ on $\UPo$, there exists
a holomorphic line bundle $\mathfrak{L} \to \CP1$ together with an integrable connection $\nabla$ with regular singular points $u^*_1, u^*_2, \dots$,
such that the restriction of $\mathfrak{L}$ is an isomorphism $i: \mathfrak{L} \to \mathfrak{E}$ with
\begin{equation*}
 d\circ i\mid_{\UPo} = (i\otimes 1)\circ \nabla\mid_{\UPo} \;.
\end{equation*}
\item The holomorphic line bundle $\mathfrak{L}$ admits a global meromorphic section $\sigma$, so $\mathfrak{L}$ is meromorphically
trivial and the connection $\nabla$ coincides with a homomorphism defined by a global meromorphic Pfaffian system $(d-\theta)\rho=0$ where $\theta$
is a global meromorphic one-form on $\CP1$. Since $\theta = d\rho/\rho$ the curvature $d\theta$ vanishes.
Then, $d\theta$ can be prolonged to the whole $\CP1$.
\end{enumerate}
\end{fact}
 
\noindent
This implies the following result:
\begin{proposition}
 On the holomorphic anti-canonical line bundle $\omega^{-1}_{\CP1}$, there exists an integrable meromorphic connection \hoehe $\nabla$ with regular singular points
 $u^*_1, u^*_2, \dots, \infty$ such that the restriction of $\omega^{-1}_{\CP1}$ to $\UPo$ is isomorphic to the determinant line bundle $(\DET \mathrm{D})^* \to \UPo$ 
 of the signature operator along the fiber of $\Z \to \UPo$. The curvature $\Omega$ of $\nabla$ equals \hoehe
 \beq
  \frac{i \, \Omega}{2\pi} =  \sum_n \frac{k_n}{6} \, \delta(u-u^*_n) \; du\wedge d\bar{u} + \frac{10-N_f}{6} \; \delta(v) \; dv\wedge d\bar{v}
\eeq
whence $\int_{\CP1} i\, \Omega/(2\pi)=2$.
\end{proposition}
\begin{proof}
The determinant line bundle $(\DET \mathrm{D})^* \to \UPo$ of the signature operator along the fiber
takes the place of the holomorphic flat line bundle $\mathfrak{E}$ in Fact~\ref{roehrl}. The bundle
$\mathfrak{E}$ has the meromorphic section $\underline{\sigma}=\pmb{\Delta}^{1/6} \, (dz)^{2}$ with $\Vert \underline{\sigma}\Vert_{Q^*}=1$.
Hence, the Bismut-Freed connection acts on $\underline{\sigma}$ simply as the exterior derivative $d$. Outside the set of singular points, the determinant line bundle is isomorphic 
to $\mathfrak{L}= \omega^{-1}_{\CP1}$. In the chart $[u:1]\in \CP1$, the isomorphism is given by multiplication with $\rho_u=\pmb{\Delta}_u^{1/6}$ and identifying 
$du^{-1}$ with $(dz_u)^{2}$ such that $i (du^{-1} ) =  
\rho_u \, (dz)^{2}$. The bundle $\mathfrak{L}$ carries the integrable meromorphic connection 
\beq
 \nabla du^{-1} = - \frac{d\rho_u}{\rho_u} \otimes du^{-1} 
\eeq
with regular singular points $u_n^*$. In particular, the form $\theta_u = d\rho_u/\rho_u$ has simple poles at every regular singular point $u_n^*$ and the counterclockwise
contour integral evaluates to
\beq
 \frac{1}{2\pi i} \oint_{u_n^*} \frac{d\rho_u}{\rho_u} = \frac{k_n}{6} \;.
\eeq
Under a change of the coordinate chart from $[u:1]$ to $[1:v]$ on $\CP1$ by $u=-1/v$ the holomorphic differential transforms as $dz_u=  - v \, dz_v$. The isomorphism is given 
by multiplication with $\rho_v=\pmb{\Delta}_v^{1/6}=v^2 \, \rho_u$ and identifying $dv^{-1}$ with $(dz_v)^{2}$  
such that  $i (dv^{-1} ) =  \rho_v \, (dz_v)^{2}$. In particular, the form $\theta_v = d\rho_v/\rho_v$ has simple pole at $v=0$ and the counterclockwise
contour integral evaluates to
\beq
 \frac{1}{2\pi i} \oint_{v=0} \frac{d\rho_v}{\rho_v} = \frac{10-N_f}{6} \;.
\eeq
The connection one-forms $\theta_u$ and $\theta_v$ patch together to give a meromorphic connection on $\CP1$:  on the intersection of the two charts we have 
$dv^{-1} = u^2 \, du^{-1}$ and
\begin{equation}
\begin{split}
 \nabla dv^{-1} & =  2 u \, du\otimes du^{-1} + u^2 \, \nabla du^{-1} =  \frac{ 2 \, \rho_u \, \dfrac{du}{u} - d\rho_u}{\rho_u} \otimes \Big[ u^2 \, du^{-1} \Big]\\
& = - \frac{d\rho_v}{\rho_v} \otimes dv^{-1} \;.
\end{split}
\end{equation}
It follows that $d\circ i\mid_{\UPo} = (i\otimes 1)\circ \nabla \mid_{\UPo}$. The curvature vanishes on all open sets. 
The curvature of the line bundle extended across the singular points is given by
\begin{eqnarray}
\label{curvature2}
  \Omega  = \bar\partial \partial \log \Vert \sigma^* \Vert_{Q^*}^2 = \bar\partial \partial \log \vert \pmb{\Delta}^{-\frac{1}{6}} \vert^2 \;.
\end{eqnarray}
We obtain
\beq
  \frac{\Omega}{2\pi i} = - \sum_{n=1}^K \frac{k_n}{6} \, \delta(u-u^*_n) \; du\wedge d\bar{u} - \frac{10-N_f}{6} \; \delta(v) \; dv\wedge d\bar{v} \;.
\eeq
The equality $\int_{\CP1} i\, \Omega/(2\pi)=2$ then follows from $\sum_n k_n = N_f +2$.
\end{proof}
\begin{remark}
Under the functor of Fact \ref{roehrl} the isomorphism class of $(\omega_{\CP1}^{-1}, \nabla)$ corresponds to the monodromy representation
of $\pi_1(\UPo)$ on the flat line bundle $\DET \; \mathrm{D}$.
\end{remark}

\section{The signature of $\Z$}
\label{signature_operator}
In the case that the elliptic fibration $\Z \to \UP$ has no surface singularities and the
singular fibers $E_u = \pi^{-1}(u)$ at $u^*_1, \dots, u^*_{N_f+2}$ are nodal curves, the manifold $\Z$ is a smooth four-dimensional manifold
with boundary $\partial \Z$. The generalization of Hirzebruch's signature theorem for manifolds with boundary by Atiyah, Patodi, Singer \cite{APS} 
shows that the elliptic signature operator $\mathrm{A}$ on $\Z$ has an 
analytic index if one imposes APS boundary conditions on $\partial \Z$. The operator $\mathrm{A}$ on $\Z$ is of the form 
\beq
 \mathrm{A} = \sigma \left( \frac{\partial}{\partial |v|} + \mathfrak{D}\right)
\eeq 
near the boundary, where $|v|$ is the inward normal coordinate, $\sigma$ a certain bundle isomorphism, and
$\mathfrak{D}$ is the selfadjoint signature operator on the boundary $\partial \Z$ \cite{Atiyah2}. 
For the boundary circle $\gamma_{\infty}$ around $u=\infty$, one obtains a three-dimensional manifold $W_{\infty}=\partial\Z$ fibered over 
a circle. On the boundary component $W_\infty$ the selfadjoint signature operator $\mathfrak{D}$ on the differential
forms of even degree is
\beq
 \mathfrak{D} = *d -d*: \quad C^\infty(W_\infty)\oplus \Omega^{2}(W_\infty) \to C^\infty(W_\infty) \oplus \Omega^{2}(W_\infty) \;.
\eeq
The eigenvalues of the operator $\mathfrak{D}$ can be positive $\lambda_j$ or negative $-\mu_j$.
If we set
\beq
 \zeta_{|\mathfrak{D}|}(s) & = & \sum_j \lambda_j^{-s} + \sum_j \mu_j^{-s} \;,\qquad
 \eta_{\mathfrak{D}}(s) =  \sum_j \lambda_j^{-s} - \sum_j \mu_j^{-s} \;,\\
 \zeta_{\mathfrak{D}^2}(s) & = & \sum_j \lambda_j^{-2s} + \sum_j \mu_j^{-2s} \;,
\eeq
we obtain $\zeta_{\mathfrak{D}^2}(s)=\zeta_{|\mathfrak{D}|}(2s)$. The logarithm of the regularized determinant $\ln \det \mathfrak{D}$ should equal
\cite{Singer}
\beq
 - \frac{d}{ds} \mid_{s=0} \left( \sum_j \lambda_j^{-s} + (-1)^s \, \sum_j \mu_j^{-s} \right) \;.
\eeq
Making the choice $(-1)^s = e^{i\pi s}$ we obtain
\beq
 \ln \det \mathfrak{D} & = & - \frac{d}{ds} \mid_{s=0} \left( \frac{\zeta_{|\mathfrak{D}|}+\eta_{\mathfrak{D}}}{2}
+ e^{i\pi s}  (\zeta_{|\mathfrak{D}|}-\eta_{\mathfrak{D}} ) \right) \\
&= & - \zeta'_{|\mathfrak{D}|}(0) - \frac{i\pi}{2}  \Big(\zeta_{|\mathfrak{D}|}(0) -\eta_{\mathfrak{D}}(0) \Big) \;.
\eeq
It follows that $\ln \det |\mathfrak{D}| = - \zeta'_{|\mathfrak{D}|}(0) = - \zeta'_{\mathfrak{D}^2}(2s)/2$. Since 
$\mathfrak{D}$ is a self adjoint operator on the odd-dimensional manifold $W_\infty$ it follows that $\zeta_{|\mathfrak{D}|}(0)=0$. 
We obtain
\beq
 \frac{\det \mathfrak{D}}{\det |\mathfrak{D}|} & = & \exp{ \left( -\frac{i\pi}{2}  \eta_{\mathfrak{D}}(0) \right) } \;.
\eeq
 
\noindent
It follows:
\begin{corollary}
The elliptic surface $\Z \to \UP$ satisfies
\beq
\label{signature_with_boundary}
 \operatorname{sign} \left(\Z\right) = - \eta_{\mathfrak{D}}\left(0 \right) =-N_f\;.
\eeq
\end{corollary}
\begin{proof}
For the elliptic surface $\Z$ the canonical class is minus the fiber class. It follows that $c_1(\Z)^2=0$.
The main theorem of \cite{APS} when applied to the elliptic surface $\Z \to \UP $ with $c_1^2(\Z)=0$
yields
\beq
 \operatorname{sign} \left(\Z\right) & = & \int_{\Z} \frac{c_1^2\left(\Z\right)}{3} - \eta_{\mathfrak{D}}\left(0 \right) 
= - \eta_{\mathfrak{D}}\left(0 \right)\;.
\eeq
On the other hand, the application of Lemma~\ref{signature_coeff} yields $\operatorname{sign} \left(\Z\right)=-N_f$.
\end{proof}

\section{Conclusion and outlook}
We have shown that the Seiberg-Witten family of elliptic curves 
defines a four-dimensional, Jacobian elliptic surface $\Z \to \UP$ with boundary.
The signature of $\Z$ is the analytic index of the signature operator on $\Z$ if we impose APS boundary conditions on $\partial \Z$. 
On the other hand, we can compute the index from the logarithmic monodromy of the canonical section of the flat determinant line bundle $\DET \; \mathrm{D} \to \UP$ of the 
signature operator along the fiber of $\Z \to \UP$. The signature of $\Z$ coincides with the number $N_f$ of hypermultiplets in gauge theory.

The identification of the hypermultiplets in the $N=2$ supersymmetric low energy $SU(2)$-Yang-Mills theory 
with the zero modes of the signature operator on the Jacobian elliptic surface defined by the Seiberg-Witten curve 
is interesting in the context of string theory. String theory predicts that $N=2$ supersymmetric $SU(2)$-gauge theory in four dimensions
emerges from the compactification the type IIB string on a certain $K3$-fibration $\widetilde{X}_3 \to \CP1$.
The Calabi-Yau three-fold $\widetilde{X}_3$ is determined by the gauge bundle in the heterotic string theory, and 
in the large base limit becomes $\mathbb{C} \times \Z$. Thus, we conclude that after the compactification the hypermultiplets
must arise from the string fields on the internal manifold which are the zero modes of the signature operator.

\section*{Acknowledgments}
I would like to thank Isadore Singer and David Morrison for many helpful 
discussions and a lot of encouragement.

\end{document}